

  \hsize=30pc
  \vsize=46pc
\hoffset=.4in
{1.5}
\font\bfthree=cmbx10 scaled\magstep{2}
\font\eightrm=cmr8
\def\smalltype{\let\rm=\eightrm  \baselineskip=8pt \rm}
\baselineskip=13.5pt

\font\smc=cmcsc10
\font\math=msbm10
\def\u#1{\hbox{\math#1}}
\def\vZ{{\u Z}}

\def\nin{\noindent}
\def\cen{\centerline}
\def\a{{\alpha}}
\def\b{{\beta}}
\def\g{{\gamma}}
\def\G{{\Gamma}}

\def\k{{\kappa}}

\def\f{{\phi}}

\def\w{{\omega}}

\def\al{{\aleph}}

\def\Ext{\mathop{\rm Ext}\nolimits}

\def\cp{cyclically presented}

\def\vd{valuation domain}

\def\vds{valuation domains}

\def\Ker{\mathop{\rm Ker}}

\def\bull{{\unskip\nobreak\hfil\penalty50\hskip .001pt \hbox{}\nobreak\hfil
           \vrule height 1.2ex width 1.1ex depth -.1ex
            \parfillskip=0pt\finalhyphendemerits=0\medbreak}\rm}
\def\qed{\bull}

\def\Proposition #1 {\vskip-\parskip\indent{\smc  Proposition #1}\quad \it}
\def\Theorem #1 {\vskip-\parskip\indent{\smc Theorem #1}\quad \it}
\def\Corollary #1 {\vskip-\parskip\indent{\smc  Corollary #1}\quad \it}
\def\Lemma #1 {\vskip-\parskip\indent{\smc  Lemma #1}\quad \it}
\def\Example #1 {\vskip-\parskip\indent{\smc  Example #1}\quad}
\def\Proof.{\indent {\smc  Proof.} \rm}

\topglue .5truein

\cen  {\bfthree On a non-vanishing Ext}
\bigskip

\cen { {\smc Laszlo Fuchs}  and {\smc Saharon Shelah\footnote{$^1$}{\smalltype
  The second author was supported by the German-Israeli Foundation for
                Scientific Research and Development. Publication 
766.}}}  \bigskip

\bigskip

{\smc Abstract.} The existence of \vds\  admitting non-standard uniserial
modules  for which certain Exts do not vanish was proved in [1] under Jensen's
Diamond Principle. In this note, the same is verified using the ZFC 
axioms alone.

\vskip .4truein

In the construction of large indecomposable divisible modules over
certain \vds\ $R$, the first author used the property that $R$ satisfied
  $\Ext_R^1(Q,U) \ne 0$, where $Q$ stands for the field of quotients of
$R$ (viewed as an $R$-module) and $U$ denotes any uniserial divisible
torsion $R$-module, for instance, the module $K = Q/R$; see [1].
However, both the existence of such a \vd\ $R$ and the non-vanishing
of Ext were established only under Jensen's Diamond Principle $\diamondsuit$
  (which holds true, e.g., in G\"odel's Constructible Universe).

Our present goal is to get rid of the Diamond Principle, that is, to
verify in ZFC the existence of \vds\ $R$ that admit divisible non-standard
uniserial modules and also satisfy $\Ext_R^1(Q,U) \ne 0$ for several uniserial
divisible torsion $R$-modules $U$. (For the proof of Corollary 3, one
requires only 6 such $U$.)

\medskip

We start by recalling a few relevant definitions. By a {\it \vd}
we mean a commutative
domain $R$ with 1 in which the ideals form a chain  under inclusion.
A {\it uniserial} $R$-module $U$ is defined similarly as a module
whose submodules form a chain under inclusion. $K = Q/R$ is
  a uniserial torsion
$R$-module, it is {\it divisible} in the sense that $rK = K$  holds
for all $0 \ne r \in R$. A divisible uniserial $R$-module is called
{\it standard} if it is an epic image of the uniserial module $Q$;
  otherwise it is said to be
  {\it  non-standard}. The existence of \vds\ admitting non-standard
uniserials has been established in ZFC; see e.g. [3], [2], and the
  literature quoted there.

  As  the $R$-module $Q$ is uniserial,  it can be represented
  as the union of a well-ordered ascending chain of cyclic submodules:
\medskip \cen {$R = Rr_0 < Rr_1^{-1} < \dots < Rr_\a^{-1} < \dots
  < \bigcup_{\a < \k} Rr_\a^{-1} = Q
  \qquad (\a < \k)$,} \medskip \nin where $r_0 = 1$, $r_\a \in R$, and
$\k$ denotes an infinite cardinal and also the initial ordinal of the
same cardinality. As a consequence,  $K = \bigcup_{\a < \k}
  (Rr_\a^{-1}/R)$ where
$Rr_\a^{-1}/R \cong R/Rr_\a$ are \cp\ $R$-modules. We denote by
$\iota_\a^\b : Rr_\a^{-1}/R \to Rr_\b^{-1}/R$  the inclusion map
for $\a < \b$, and may view $K$ as the direct limit of its
submodules $Rr_\a^{-1}/R$ with the monomorphisms $\iota_\a^\b $
  as connecting maps.

A uniserial divisible torsion module $U$ is a {\it `clone'} of
  $K$ in the  sense of Fuchs-Salce [2], if there are units
$e_\a^\b \in R$ for all pairs $\a < \b \ (< \k)$ such that
\medskip \cen {$e_\a^\b e_\b^\g - e_\a^\g \in Rr_\a$
\quad for all $\a < \b < \g < \k,$} \medskip \nin and $U$ is the
  direct limit of the direct system
of the modules $Rr_\a^{-1}/R$ with connecting maps $ \iota_\a^\b e_\a^\b :
Rr_\a^{-1}/R \to Rr_\b^{-1}/R \ (\a < \b) $; i.e. multiplication by
$e_\a^\b$ followed by the inclusion map. It might be helpful to point
out that though $K$ and $U$ need not be isomorphic, they are `piecewise'
isomorphic in the sense that they are unions of isomorphic pieces.

Let $R$ denote the \vd\ constructed in the paper [1] (see also
Fuchs-Salce [2]) that satisfies $\Ext_R^1(Q,K) \ne 0$ in the constructible
  universe L. Moreover, there are clones $U_n$ of $K$, for any integer $n>0$,
  that satisfy $\Ext_R^1(Q,U_n) \ne 0$: for convenience, we let $K=U_0$.

This $R$ has the value group  $\G = \oplus_{\a < \w_1} \vZ$,
ordered anti-lexicographically, and its quotient field $Q$ consists of all
  formal rational functions of $u^\g$ with coefficients in an arbitrarily
chosen, but fixed field, where
$u$ is an indeterminate and $\g \in \G$. It is shown in [2] that such
an $R$ admits divisible non-standard uniserials (i.e.
clones of $K$ non-isomorphic to $K$), and under the additional hypothesis of
  $\diamondsuit_{\al_1}$, $\Ext_R^1(Q,U_n) \ne 0$ holds; in other 
words, there is a
non-splitting exact sequence \medskip \cen
{$0 \to U_n \to H_n\ {\buildrel \f \over \longrightarrow}\ Q \to 0$.} \medskip

Using the elements $r_\a \in R$ introduced above, for each $n < \w$ we
  define a tree $T_n$ of length $\k$  whose set of vertices at level
  $\a$ is given by
\medskip \cen {$ T_{n\a} = \{x \in H_n \ | \ \f (x) = r_\a^{-1} \}.$ }
   \medskip \nin
The partial order $<_T$ is defined in the following way: $x <_T y$ in $T_n$
   if and only if, for some $\a < \b$, we have $\f x = r_\a^{-1} $ and
$\f y = r_\b^{-1} $ such that
\medskip \cen {$ x = r_\a^{-1}  r_\b  y$ \quad  in $H_n$,}\medskip \nin
where evidently $r_\a^{-1}  r_\b  \in R$.

Fix an integer $n >0$, and define $T$ as the union of the trees
$T_0, T_1, \dots, T_n$ with a minimum element $z$ adjoined. It is
straightforward to check that
  $(T, <_T )$ is indeed a tree with $\k$ levels, and the
inequalities $\Ext_R^1(Q,U_i) \ne 0 \ (i=0,1, \dots, n)$ guarantee that
$T$ has no branch of length $\k$.

We now define a model {\bf M} as follows. Its universe is the union of the
  universes of  $R, Q, U_i, H_i \ (i=0, \dots, n),$ and it has the following
  relations:

(i) unary relations $R,Q,U_i,H_i,T$, and   $S = \{r_\a \ | \ \a<\k \}$,

(ii) binary relation $<_T$, and $<_S$ (which is the natural 
well-ordering on $S$),

(iii) individual constants $0_R, 0_Q, 0_{U_i}, 1_R$, and

(iv) functions: the operations in $R, Q, U_i, H_i$, where $R$ is a domain,
$Q, U_i, H_i$ are  $R$-modules, $\f_i$ is an $R$-homomorphism from $H_i$ onto
$Q \ (i=0, \dots, n$,
and $\psi: Q \to K$ with $\Ker \psi = R$ is the canonical map.

  We argue that even though our universe V does not satisfy V=L, the class L
does satisfy it, and so in L we can define
the model {\bf M} as above.
Let {\bf T} be the first order theory of {\bf M}.
So the first order (countable theory) {\bf T} has in L a model in which

$(^\star)_{\bf M}$  the tree $(T,<_T)$
with set of levels $(S,<_S)$ and with  the function $\f = \cup \f_i)$
  giving the levels, as
interpreted in {\bf M}, has no full branch (this means that there is 
no function
  from $S$ to $T$ increasing in the natural sense and inverting $\f$,
ot any $\f_i$).

Hence we conclude as in [3] (by making use of Shelah [4]) that
there is a model ${\bf M'}$
  with those properties in V (in fact, one of cardinality $\al_1$).

Note that $(S,<_S)$ as interpreted in ${\bf M'}$ is not well ordered,
but  it is still a linear order of uncountable cofinality (in fact, of
cofinality  $\al_1$),
  the property $(^\star)_{\bf M'}$  still holds,
and it is a model of {\bf T}.
  This shows that all relevant properties of {\bf M} in L hold for
${\bf M'}$ in V, just as indicated in [3].

It should be pointed out that, as an alternative, instead of using a
smaller universe of set theory L,
we could use a generic extension not adding new subsets of the natural numbers
(hence essentially not adding new countable first order theories like
{\bf T}).

If we continue with the same argument as in [3], then using [4] we can
  claim that we have proved in ZFC the following theorem:

\Theorem 1. There exist \vds\ $R$ admitting non-standard uniserial
torsion divisible modules
such that $\Ext_R^1(Q,U_i) \ne 0$ for various clones $U_i$ of $K$.  \qed

Hence we derive at once that the following two corollaries
are true statements in ZFC; for their proofs we refer to [2].

\Corollary 2. There exist \vds\ $R$ such that if $U, V$ are non-isomorphic
clones of $K$, then $\Ext_R^1(U,V)$
satisfies:

{\rm (i)} it is a divisible mixed $R$-module;

  {\rm (ii)} its torsion submodule is uniserial.   \qed

More relevant consequences are stated in the following corollaries; they solve
Problem 27 stated in [2].

\Corollary 3. There exist \vds\  admitting indecomposable divisible modules
of cardinality larger than any prescribed cardinal.  \qed

\Corollary 4. There exist \vds\ with superdecomposable divisible modules
of countable Goldie dimension.   \qed

\bigskip \bigskip

\cen {\bf References} \medskip

[1] {\smc L. Fuchs}, Arbitrarily large indecomposable divisible
  torsion modules over certain
valuation domains, {\sl Rend. Sem. Mat. Univ. Padova} 76 (1986), 247-254.

[2] {\smc L. Fuchs} and {\smc L. Salce}, Modules over non-Noetherian Domains,
{\sl Math. Surveys and Monographs,} vol. 84, American Math. Society
(Providence, R.I., 2001).

[3] {\smc L. Fuchs} and {\smc S. Shelah},
  Kaplansky's problem on valuation rings, {\sl Proc. Amer. Math. Soc.}
  105 (1989), 25-30.

[4] {\smc S. Shelah}, Models with second order properties. II. Trees with
no undefined branches, {\sl Ann. of Math. Logic} 14 (1978), 73-87.

\vskip .5 truein

{\smc \nin Department of Mathematics

\nin Tulane University,
New Orleans, Louisiana 70118, USA}

\nin e-mail: fuchs@tulane.edu

\medskip and \medskip

{\smc \nin Department of Mathematics

\nin Hebrew University, Jerusalem, Israel 91904

\nin {\rm and} Rutgers University, New Brunswick, New Jersey 08903, USA}

\nin e-mail: shelah@math.huji.ac.il
\bye